\documentstyle[11pt]{article}

\newtheorem{thm}{Theorem}[section]
\newtheorem{lem}[thm]{Lemma}

\newtheorem{pro}[thm]{Proposition}
\newtheorem{ex}[thm]{Example}

\newtheorem{defi}[thm]{Definition}

\newcommand{\gm }{\Gamma }
\newcommand{\lon }{\longrightarrow }
\newcommand{\be }{\begin{eqnarray*}}
\newcommand{\ee }{\end{eqnarray*}}

\newcommand{\per }{\backl }

\input amssym.def 
\input amssym

\newcommand{\pf}{\noindent{\bf Proof.}\ }
\newcommand{\qed}{\begin{flushright} $\Box$\ \ \ \ \ \
                  \end{flushright}}

\newcommand{\reals}{{\Bbb R}}

\newcommand{\frakg}{{\frak g}}

\newcommand{\integers}{{\Bbb Z}}

\newcommand{\backl}{\mathbin{\vrule width1.5ex height.4pt\vrule height1.5ex}}
\newcommand{\cala}{{\cal A}}

\newcommand{\smalcirc}{\mbox{\tiny{$\circ $}}}

\def\description label#1{\hfil\bf[#1]\hfil}
\newcommand{\secc}[1]{\gm(\wedge^{#1}A)}
\newcommand{\secd}[1]{\gm(\wedge^{#1}A^{*})}
\newcommand{\gma}[1]{\gm (\wedge^{ #1}A)}
\newcommand{\gmd}[1]{\gm (\wedge^{ #1}A^{*} )}
\newcommand{\Ger}{Gerstenhaber  }
\newcommand{\Koszul}{generating  }
\newcommand{\equivalent}{homotopic }

\def\sdp{\mathbin{\hbox{$\mapstochar\kern-.3333em\times$}}}
\def\pds{\mathbin{\hbox{$\times\kern-.55em\mapstochar\,$}}}

\newcommand{\wed}{\mathbin{\lower1.5pt\hbox{$\scriptstyle{\wedge}$}}}

\let\Tilde=\widetilde

\def\chigh{{\raise1.5pt\hbox{$\chi$}}}
\let\phi=\varphi
\def\til0{\Tilde{0}}

\def\dminus{\raise2pt\hbox{\vrule height1pt width 2ex}\hskip3pt}

\def\pback#1{\mathbin{{{\lower1.2ex\hbox{$\times$}}\atop #1}}}

\def\vlra{\hbox{$\,-\!\!\!-\!\!\!-\!\!\!-\!\!\!-\!\!\!
-\!\!\!-\!\!\!-\!\!\!-\!\!\!-\!\!\!\longrightarrow\,$}}

\def\gpd{\,\lower1pt\hbox{$\longrightarrow$}\hskip-.24in\raise2pt
             \hbox{$\longrightarrow$}\,}

\def\lgpd{\,\lower1pt\hbox{$\vlra$}\hskip-1.02in\raise2pt\hbox{$\vlra$}\,}

\def\llgpd{\,\lower1pt\hbox{$\vvlra$}\hskip-1.3in\raise2pt\hbox{$\vvlra$}\,}

\begin{document}

\title{\bf \Ger algebras and  BV-algebras in  Poisson geometry}

\author{ PING XU \thanks{ Research partially supported by NSF
        grant DMS95-04913.}\\
 Department of Mathematics\\The  Pennsylvania State University \\
University Park, PA 16802, USA\\
        {\sf email: ping@math.psu.edu }}

\date{February, 1997}

\maketitle

\begin{abstract}
The purpose of this paper is to establish an explicit correspondence
between various geometric structures on a 
vector bundle with some well-known algebraic  structures
such as \Ger algebras and BV-algebras. Some applications
are discussed. In particular, we found an explicit 
connection  between the  Koszul-Brylinski operator of a Poisson manifold
and its  modular class. As a consequence, we   prove that Poisson homology is
isomorphic to Poisson cohomology for
unimodular Poisson structures.
\end{abstract}

\section{Introduction}

BV-algebras arise from the BRST theory of topological
field theory \cite{Witten}. 
Recently, there has been  a great deal of interest  in these algebras
in connection with various  subjects such as operads  and string theory
 \cite{G} \cite{GK} \cite{LZ}  \cite{KVS} \cite{PS} \cite{S1} \cite{S2}
 \cite{Z}.

Let us first  recall various relevant definitions below.  
Here we will follow those as  in \cite{KS0}.

 A {\em  Gerstenhaber algebra} consists of a  triple
$(\cala =\oplus_{i\in \integers}\cala^{i}, \wedge, [\cdot , \cdot ] )$ such that
$(\cala , \wedge )$ is a graded commutative associative algebra,
and $(\cala=\oplus_{i\in \integers}\cala^{(i)}, [\cdot , \cdot ])$,
with $\cala^{(i)}=\cala^{i+1}$,  is  a  graded Lie algebra, and
 $[a, \cdot ]$, for each $a\in  \cala^{(i)}$ is a derivation  with
respect to $\wedge$
with degree $i$.

An operator $D$ of degree $-1$ is said to generate the Gerstenhaber
algebra bracket if for every $a\in \cala^{|a|}$ and
$b\in \cala$,

\begin{equation}
[a, b]=(-1)^{|a|}(D(a\wedge b)-Da\wedge b -(-1)^{|a|}a\wedge Db) .
\end{equation}

A Gerstenhaber algebra is  said to be {\em exact} if there is
an operator $D$ of square zero generating 
the bracket. In this case,  $D$ is called a {\em \Koszul operator}.
An  exact Gerstenhaber algebra is also called a {\em Batalin-Vilkovisky algebra}
(or {\em BV-algebra} in short).

A {\em differential Gerstenhaber algebra} is a Gerstenhaber algebra
equipped with a differential $d$, which is a derivation
of degree $1$  with respect to $\wedge$ and $d^2=0$.
 It is called a {\em strong
differential \Ger algebra} if, in addition,  $d$ is
derivation of the graded Lie bracket.


Kosmann-Schwarzbach noted \cite{KS}
that  these algebra structures had also appeared
in Koszul's work in 1985 \cite{Koszul} in his study of Poisson  manifolds.
In fact, as pointed out in \cite{KS}, these  examples 
are connected with a certain  differential  structure on vector bundles, 
called {\em Lie algebroids} by  Pradines \cite{Pradines2}. Let us 
 recall for the
benefit of the reader the definition of a Lie algebroid \cite{Pradines2} 
\cite{Pradines1}.
 
\begin{defi}\label{Lie.Algebroid} A  Lie algebroid $A$
over a manifold $M$ is a vector bundle $A$ over $M$ together with  a Lie
algebra structure on the space $\Gamma(A)$ of smooth sections of $A$,
and  a bundle map $a: A \rightarrow TP$ (called the anchor),
 extended to a map between sections of these bundles,  such that
 
(i) $a([X,Y])=[a(X),\  a(Y)]$; and
 
(ii) $[X, fY] = f[X,Y] + (a(X) f)Y$
 
 for any smooth sections $X$ and $Y$ of $A$ and any smooth function
$f$ on $M$.
\end{defi}

Among many examples of Lie algebroids are
 usual Lie algebras, the tangent bundle
of a manifold, and an integrable distribution over a manifold
 (see \cite{Mackenzie}). 
Recently,    Lie algebroids attract  increasing interest in  Poisson
geometry. One of the main reason is due to the following example
connected with Poisson manifolds.  

Let $P$ be  a  Poisson manifold with Poisson tensor
  $\pi $. Then  $T^{*}P$  inherits   a  natural
Lie algebroid structure, called
the cotangent Lie algebroid of the Poisson manifold $P$ \cite{CDW}.
 The anchor map
 $\pi^{\#}: T^*P \rightarrow TP$ is defined by
\begin{equation}
\label{eq:anchor0} 
 \pi^{\#}: ~~T_{p}^{*}P \longrightarrow T_p P:~~
\pi^{\#}(\xi)(\eta )=\pi (\xi ,\eta ), \ \ \ \forall \xi , \eta\in T_{p}^{*}P
\end{equation}
 and the Lie bracket of $1$-forms $\alpha$ and $\beta$ is given by
 \begin{eqnarray}
 \label{eq_bracket-on-one-forms}
 [\alpha, \beta  ]
 & = & - d \pi (\alpha, \beta) ~ + ~ L_{{\pi}^{\#}(\alpha)} \beta ~ -
 ~ L_{{\pi}^{\#}(\beta)} \alpha.
 \end{eqnarray}

In  \cite{KS},  Kosmann-Schwarzbach constructed various
examples of strong differential  \Ger algebras and  BV-algebras
in connection with  Lie algebroids.
Motivated by \cite{KS}, in this note, we will establish
a more precise relation between these
algebra structures and  some of the well-known
geometric structures in Poisson geometry.

More precisely, we will investigate the following question:
Let $A$ be a  vector bundle of rank $n$   over the base $M$,
and let $\cala =\oplus_{0\leq k\leq n} \gma{k}$.  With
respect to the wedge product, $\cala$ becomes a  graded commutative
associative algebra. Then  the question is: 

\begin{quote}
What additional structure on $A$ will make $\cala$
into  a \Ger algebra,  a strong differential \Ger
algebra, or  an  exact \Ger algebra (or a  BV-algebra)?
\end{quote}

The answer is surprisingly  simple. 
\Ger algebras  and  strong differential  \Ger algebras,
correspond,  exactly to the structures of
  Lie algebroids  and  Lie bialgebroids (see Section 2
for the definition), respectively,
 as already  indicated in \cite{KS}.
And  an  exact \Ger algebra  structure
corresponds to a  Lie  algebroid $A$  together
 with a flat $A$-connection on its canonical line bundle
$\wedge^{n}A$.  
This fact was  implicitly
contained in Koszul's work \cite{Koszul}
although he treated only the case of multivector
fields. However, the  formulas
(\ref{eq:D}) and  (\ref{eq:connection}) establishing
the  explicit correspondence  seem to be new. 

Below is a table of  the  correspondence.

\bigskip
 
\begin{tabular}{ccc}
 
Structures on algebra $\cala$ & & Structures on the vector bundle $A$ \\ \hline \hline
 & & \\
\Ger algebras  & $\leftrightarrow $ & Lie algebroids \\
 & &  \\
strong differential \Ger algebras & $\leftrightarrow $ &
Lie bialgebroids \\
 & &  \\
exact \Ger algebras (BV-algebra) & $\leftrightarrow $ &  Lie algebroids with a flat
                                           $A$-connection on $\wedge^{n}A$ 
        \end{tabular}
 
\bigskip

The content above occupies Section 2 and Section 3. Section 4
is devoted to applications. In particular, 
we establish an explicit connection between  the  Koszul-Brylinski operator
 on a Poisson manifold  with its  modular class.
As a consequence, we   prove that Poisson homology is
isomorphic to Poisson cohomology for
unimodular Poisson structures (see \cite{we:modular}  \cite{BZ} 
for the definition).

As  another application, we introduce the notion of
 Lie algebroid homologies, which
are the homology groups induced by \Koszul operators 
$D: \secc{*}\lon \secc{*-1}$.
Since a  \Koszul operator on a Lie algebroid depends on
the choice of a flat $A$-connection on the canonical line
bundle $\wedge^{n}A$, in general the homology depends on the
choice of such a connection $\nabla$. When two connections
are \equivalent (see Section 4 for the precise definition),
 their corresponding  homology groups are isomorphic.
So for a given Lie algebroid, its homologies are    in fact parameterized
by the first Lie algebroid cohomology $H^{1}(A, \reals )$.
When $A$ is a Lie algebra and $\nabla$ is the trivial connection,
this reduces to the usual Lie algebra homology with trivial
coefficients. On the other hand, Poisson homology
can also be considered as a special case of  Lie algebroid homology,
when $A$ is taken as  the cotangent Lie algebroid of a Poisson
manifold.

We note that in a recent preprint \cite{EvansLW}, Evens, Lu and Weinstein have
also established a connection between 
Poisson homology and the  modular class of Poisson manifolds.
Finally we also would like to refer the reader to   a recent preprint
\cite{H}  of   Huebschmann for its  close connection with the present paper. 
 
{\bf Acknowledgments.}  The author  would like to thank 
 Jean-Luc Brylinski,  Jiang-hua Lu, and Alan Weinstein for useful discussions.
He is especially grateful to Yvette Kosmann-Schawrzbach and Jim
Stasheff for providing
 many helpful comments on the first draft of the manuscript.
In addition to the funding sources mentioned
in the first footnote, he would also  like to thank 
IHES  and
  Max-Planck Institut for their hospitality and financial support while part of this project was being done.

\section{Gerstenhaber algebras and differential Gerstenhaber algebras}

In this section, we will treat \Ger algebras and differential
\Ger algebras arising from a vector bundle.

Again, let $A$ be a  vector bundle of rank $n$ over  $M$,
and let $\cala =\oplus_{0\leq k\leq n} \gma{k}$.  The following
proposition establishes a one-one correspondence between \Ger algebra
structures on $\cala$ and Lie algebroid structures on  the
vector bundle $A$.

\begin{pro}
\label{pro:ger}
 $\cala$ is a Gerstenhaber algebra iff $A$ is a Lie algebroid.
\end{pro}

This is a well-known result (see \cite{GS} \cite{KS-M} \cite{MX}).
 For completeness, we will sketch a proof
below.\\\\
\pf   Suppose that there is a graded Lie bracket $[\cdot , \cdot ]$
that makes $\cala$ into a \Ger algebra. 
It is clear that $(\gm (A), [\cdot , \cdot ])$  is   a Lie
algebra.

Second, for any $X\in \gm (A)$ and $f, g\in C^{\infty}(M)$,
it follows from the derivation property that
$$[X, fg]=[X, f]g+f[X, g]. $$
Hence, $[X, \cdot ]$ defines a vector field on $M$,
which will be denoted by $a (X)$.
It is easy to see that $a$ is  in fact   induced by
a bundle map from $A$ to $TP$. By applying the  graded Jacobi identity,
it follows that 
$$a ([X, Y])=[a  (X) ,  a (Y) ]. $$

Finally, again from the derivation property, it follows
that
$$[X, fY]=(a (X)f)Y+f [X, Y] .$$ 

This shows that $A$ is indeed  a Lie algebroid.

Conversely, given a Lie algebroid $A$, it is easy to check  that
$\cala =\oplus_{0\leq k\leq n} \gma{k}$  forms a \Ger algebra
(see \cite{KS} \cite{MX}). \qed

The following lemma   gives an alternative  way of
characterizing a Lie  algebroid, which  
should  be of interest itself.

Recall that a {\em differential graded algebra}
is a graded commutative associative algebra equipped
with  a differential $d$, which is a 
derivation of degree $1$ and of square zero.

\begin{lem} 
\label{lem:diff}
\cite{KS-M} \cite{Ko}
Given a vector bundle $A$ over  $M$. 
$A$ is a Lie algebroid iff $\wedge^{*} A^*$ is
a differential graded algebra.
\end{lem}
\pf  Given a  Lie algebroid $A$, it is known that
the space of sections $\gm (\wedge^{*} A^* )$   admits  a 
differential  $d$ that  makes it into a 
differential graded algebra \cite{KS-M}. In this case, $d: \gm (\wedge^{k} A^* )
\lon \gm (\wedge^{k+1} A^* )$ is  simply the differential defining
the Lie algebroid cohomology given  as  below (see \cite{Mackenzie}
\cite{MX}  \cite{WX}):
\begin{eqnarray}
d\omega (X_1,\ldots ,X_{k+1}) & = & \sum_{i=1}^{k+1}  (-1)^{i+1} a(X_i)
(\omega (X_1,\hat{\ldots}, X_{k+1})) \nonumber \\
& & \qquad + \sum_{i\ < j}  (-1)^{i+j} \omega ([X_i ,X_j ],X_1 ,
\hat{\ldots}\,\hat{\ldots}\,, X_{k+1}), \label{eq:derivative}
\end{eqnarray}
for $\omega  \in \secd{k}$, $X_i \in \Gamma A,\ 1\leq i\leq k+1$.

Conversely, if $\gm (\wedge^{*} A^* )$ is a differential
graded algebra  with  differential $d$,
then the equations:
\begin{equation}
\label{eq:anchor}
a (X)f=<df , X>, \ \ \ \forall f\in C^{\infty}(M) \ \mbox{\ and }
X\in  \gm (A), 
\end{equation}
and
\begin{equation}
\label{eq:bracket0}
<[X, Y], \theta >=a(X)(\theta \cdot Y)- a(Y)(\theta \cdot X)
-(d\theta ) (X, Y) 
\end{equation}
define the anchor map and the Lie bracket of
a Lie algebroid structure on $A$. \qed
{\bf Remark.} The lemma above is essentially Proposition 6.1
of \cite{KS-M}. Equation (\ref{eq:bracket0}) is Formula (6.6) in
\cite{KS-M}. \\\\

Recall that a {\em Lie bialgebroid} \cite{KS}\cite{MX} is
 a dual pair $(A,A^*)$ of vector bundles 
equipped with Lie algebroid structures such that the differential $d_*$, 
induced  from the  Lie algebroid  structure on $A^*$ as defined 
by Equation (\ref{eq:derivative}),  is a derivation 
of the Lie bracket on    $\gm (A )$, i.e.,
\begin{equation}
d_{*}[X, Y]=[d_{*}X, Y]+[X, d_{*}Y], \ \ \forall X ,Y \in \gm (A).
\end{equation}

The following result is due to  Kosmann-Schwarzbach \cite{KS}. 
\begin{pro}
\label{pro:dG}
$\cala$ is a  strong differential
Gerstenhaber algebra iff $A$ is a Lie bialgebroid.
\end{pro}
\pf Assume that $\cala$ is a strong differential \Ger algebra.
Then, $A^*$ is a Lie algebroid according to Lemma \ref{lem:diff}.
Moreover, the
derivation property  of the differential 
with respect to  the Lie bracket on $\gm (A )$ 
implies  that $(A^{*}, A)$ is a Lie bialgebroid. 
This is equivalent to that $(A, A^{*})$ is a Lie bialgebroid
by duality \cite{MX}.
Conversely, it is straightforward
to see, for a given  Lie bialgebroid $(A, A^{*})$,
that $\cala$ is a strong differential \Ger algebra (see \cite{KS}). \qed

\begin{ex}
Let $P$ be a Poisson manifold with Poisson tensor 
$\pi$.  Let $A=TP$ with the tangent algebroid structure.
 It is well known that the space of
multivector fields $\cala =\oplus (\wedge^{k}TP)$
has a Gerstenhaber algebra structure, where the
graded Lie bracket is called  the   Schouten bracket.

In 1977, Lichnerowicz introduced a differential
operator $d=[\pi , \cdot ]$, which he used to define
the Poisson cohomology \cite{Lich}. 
It is obvious that $\cala$ becomes a  strong
differential Gerstenhaber algebra, so it
should correspond  to a Lie bialgebroid structure on
$(TP, T^{*}P)$ according to  Proposition \ref{pro:dG}. It is not surprising
that this is just  the standard Lie
bialgebroid of a Poisson manifold \cite{MX}, where
the Lie algebroid structure on $T^{*}P$ is defined as in  the introduction
(see Equations (\ref{eq:anchor0}) and (\ref{eq_bracket-on-one-forms})).
It is, however, quite  amazing that the Lie algebroid structure
on $T^*P$ was not known until the middle
of 1980's  (see \cite{KS1} for the references)
and the Lie bialgebroid structure comes  much later! For the Lie algebroid
$T^{*}P$, the associated differential operator on $\gm(\wedge^{*}TP)$ is the
Lichnerowicz differential $d=[\pi , \cdot ]$. This property was proved,
independently by Bhaskara and Viswanath \cite{BV}, and Kosmann-Schwarzbach
and Magri \cite{KS-M}.
\end{ex} 

\section{Exact Gerstenhaber algebra structures  on the exterior
algebra of a vector bundle}

In this section, we will  move to exact \Ger algebras arising from 
a vector bundle.

Let $A\lon M$ be a  Lie algebroid with anchor $a$ and $E\lon M$
a vector bundle over $M$. By an $A$-connection on $E$, we
mean a $\reals$-linear  map: 
$$\gm (A)\otimes \gm (E)\lon \gm (E)$$
$$X\otimes s\lon \nabla_{X}s$$
satisfying the axioms resembling   those of the usual  
 linear connections, i.e.,
$\forall f\in C^{\infty}(M), \ X\in \gm (A), s\in \gm (E)$,
\be
&&\nabla_{fX}s=f\nabla_{X}s;\\
&&\nabla_{X}(fs)=(a(X)f) s+f \nabla_{X}s .
\ee
Similarly, the curvature $R$ of an $A$-connection $\nabla$
is the element in $\gm (\wedge^{2}A^{*})\otimes End (E)$
defined by
\begin{equation}
R(X, Y)=\nabla_{X}\nabla_{Y}- \nabla_{Y}\nabla_{X}-\nabla_{[X, Y]}, \ \ 
\forall  X, Y\in \gm (A). 
\end{equation}

Given a  Lie algebroid $A$ and an $A$-connection  $\nabla$ on the canonical line 
bundle $E=\wedge^{n}A$, 
define a differential operator $D: \gm (\wedge^{k}A) \lon
 \gm (\wedge^{k-1}A)$
as follows. Let $U$ be any section in $\gm (\wedge^{k}A)$ and
write, locally,
 $U=\omega \per \Lambda$, where $\omega \in \gm (\wedge^{n-k}A^{*})$
and $\Lambda \in \gm (\wedge^{n}A )$.
At each point $m\in M$, set 
\begin{equation}
\label{eq:D}
DU|_{m}=- (-1)^{|\omega |}( d\omega 
\per \Lambda+\sum_{i=1}^{n}(\alpha_{i}\wedge \omega )
\per \nabla_{X_{i}}\Lambda ) ,
\end{equation}
where $X_{1}, \cdots , X_{n}$ is a basis of $A|_{m}$ and
$\alpha_{1}, \cdots , \alpha_{n}$ its dual basis
in $A^{*}|_{m}$. Clearly,
this definition is independent of the choice of the basis. \\\\
{\bf Remark} We would like to make a remark on the notation.
Let $E$ be a vector bundle over $M$. Assume that $V\in \gm  (\wedge^{k}E)$
and $\theta \in (\wedge^{l}E^{*})$ with $k\geq l$. Then,
by $\theta \per V$  we denote the  section of $\wedge^{k-l}E$
given by 
$$(\theta \per V)(\omega )=V(\theta\wedge \omega ), \ \ \ , \forall \omega
\in \gm (\wedge^{k-l}E^{*}).$$
We will stick to this notation in the sequel no matter
whether  $E$ is the Lie algebroid $A$ itself or its
dual $A^{*}$.

\begin{pro}
\label{pro:D2}
$D$ is a well-defined operator and 
$$D^{2}U=-R\per U, $$
where $R\in \gm (\wedge^{2}A^{*})$ is  the curvature
of the connection $\nabla$ (note that $End E$ is a trivial
line bundle).
\end{pro}
\pf  Suppose that $f$ is a locally
 nonzero function on $M$,
and $U=f\omega \per \frac{1}{f} \Lambda $.
Then,
\be
&& d(f\omega )\per \frac{1}{f}  \Lambda + \sum_{i}(\alpha_{i}\wedge f\omega )  
\per \nabla_{X_{i}}( \frac{1}{f}\Lambda )\\
&=&\frac{1}{f} (df \wedge \omega ) \per \Lambda  +d\omega \per \Lambda\\ 
&&+\sum_{i} f ((aX_{i})\frac{1}{f}  )(\alpha_{i}\wedge \omega ) \per \Lambda
+\sum_{i} (\alpha_{i}\wedge \omega  )\per  \nabla_{X_{i}} \Lambda \\
&=&\frac{1}{f} (df \wedge \omega ) \per \Lambda   +
f(d(\frac{1}{f}) \wedge \omega ) \per \Lambda + d\omega \per \Lambda
+\sum_{i} (\alpha_{i}\wedge \omega ) \per  \nabla_{X_{i}} \Lambda\\
&=& d\omega \per \Lambda
+\sum_{i} ( \alpha_{i}\wedge \omega ) \per  \nabla_{X_{i}} \Lambda.
\ee
Therefore, $D$ is well-defined.

For the second part, we have

\be
D^{2}U&=&-(-1)^{|\omega |}D (
d\omega \per \Lambda+\sum_{i=1}^{n}(\alpha_{i}\wedge \omega )
\per \nabla_{X_{i}}\Lambda )\\
&=&-(\sum_{i}(\alpha_{i}\wedge d\omega ) \per \nabla_{X_{i}}\Lambda
+\sum_{i}d(\alpha_{i}\wedge \omega )\per \nabla_{X_{i}}\Lambda
+\sum_{j, i}(\alpha_{j}\wedge \alpha_{i}\wedge \omega )
\per \nabla_{X_{i}}\nabla_{X_{j}}\Lambda  )\\
&=& -(\sum_{i}(d\alpha_{i}\wedge \omega ) \per \nabla_{X_{i}}\Lambda
+\sum_{j, i} (\alpha_{j}\wedge \alpha_{i}\wedge \omega )
\per \nabla_{X_{j}}\nabla_{X_{i}}\Lambda  )\\
&=& -[ \sum_{i} \omega \per (d\alpha_{i} \per \nabla_{X_{i}}\Lambda )
+\sum_{j, i} \omega \per (\alpha_{i}\wedge \alpha_{j}
\per \nabla_{X_{i}}\nabla_{X_{j}}\Lambda  )]
\ee

The conclusion thus  follows from the following

\begin{lem}
$$ \sum_{i} d\alpha_{i} \per \nabla_{X_{i}}\Lambda 
+\sum_{j, i} (\alpha_{i}\wedge \alpha_{j} )
\per \nabla_{X_{i}}\nabla_{X_{j}}\Lambda =-R\per \Lambda .  $$
\end{lem}
\pf  It is a straightforward verification, and is left to the
readers. \qed

\begin{pro}
\label{pro:bv}
Let  $D: \gm (\wedge^{k}A) \lon \gm (\wedge^{k-1}A)$ be the operator
 defined in Equation (\ref{eq:D}). Then,
$D$ generates the Gerstenhaber algebra  bracket on  $\gm (\wedge^{*}A )$,
i.e, for any $U\in \gm (\wedge^{u}A) $ and $V\in \gma{v}$,
\begin{equation}
[U, V]=(-1)^{u}(D(U\wedge V)-DU\wedge V- (-1)^{u}U\wedge DV ).
\end{equation}
\end{pro}

We need  a couple of  lemmas before  proving  this proposition.

\begin{lem}
For any $U\in \gma{u}, \ V\in \gma{v}$ and $\theta \in \gmd{u+v-1}$,
\begin{equation}
\label{eq:bracket}
[U, V]\per \theta =(-1)^{(u-1)(v-1)}U\per d(V\per \theta )-
V\per d(U\per \theta ) -(-1)^{u+1}(U\wedge V)\per d\theta  .
\end{equation}
\end{lem}
\pf See Equation (1.16) in  \cite{Vaisman}.

\begin{lem}
For any $U\in \gma{u}$ and $\theta \in \gmd{u-1}$,
\begin{equation}
\label{eq:contraction}
\theta \per DU=(-1)^{|\theta |}D(\theta \per U)+d\theta \per U .
\end{equation}
\end{lem}
\pf Assume that $U=\omega \per \Lambda $. Then, $\theta \per
U=(\omega \wedge \theta ) \per \Lambda $, and therefore,

\be
D(\theta \per U)&=&-(-1)^{|\omega |+|\theta |}(d (\omega \wedge \theta )\per
\Lambda +\sum_{i}(\alpha_{i} \wedge \omega  \wedge \theta ) \per
\nabla_{X_{i}}\Lambda )\\
&=& -(-1)^{|\omega |+|\theta |} (d \omega \wedge \theta \per \Lambda 
+(-1)^{|\omega |}(\omega \wedge d\theta ) \per \Lambda 
+\sum_{i}(\alpha_{i}\wedge \omega \wedge \theta ) \per \nabla_{X_{i}}\Lambda )\\
&=&(-1)^{|\theta |}(\theta \per DU)-(-1)^{|\theta |}d \theta \per U.
\ee
\qed
{\bf Proof of Proposition \ref{pro:bv}} For any  $U\in \gma{u}$,
$V\in \gma{v}$ and $\theta \in \gmd{u+v-1}$,
using Equation (\ref{eq:contraction}), we have

$$\theta \per D(U\wedge V)=(-1)^{|\theta |}D(\theta \per (U\wedge V))
+d\theta \per (U\wedge V).$$

On the other hand, we have
\be
\theta \per (U\wedge DV)&=&(U\per \theta )\per DV\\
&=&(-1)^{|\theta |-u}D((U\per \theta )\per V)+d (U\per \theta )\per V,
\ee
and

\be
&&\theta \per (DU\wedge V) \\
&=&(-1)^{(u-1)v}\theta \per (V\wedge DU)\\
&=&(-1)^{(u-1)v}((-1)^{|\theta |-v}D((V\per \theta )\per U)+d (V\per \theta )\per U )\\
&=&(-1)^{uv+|\theta |}D((V\per \theta )\per U)
+(-1)^{(u-1)v}d (V\per \theta )\per U ) .
\ee

It thus follows that
\be
&&\theta \per ((-1)^{u} D(U\wedge V)-(-1)^{u}DU\wedge V-U\wedge DV)\\
&=&(-1)^{|\theta |+u}D(\theta \per (U\wedge V))-
(-1)^{uv+|\theta |+u} D((V\per \theta )\per U)
-(-1)^{|\theta |-u}D((U\per \theta )\per V)\\
&&+(-1)^{u}d\theta \per (U\wedge V) +(-1)^{(u-1)(v-1)}d(V\per \theta )\per
U-d (U\per \theta )\per V .
\ee
The conclusion thus follows from Equation  (\ref{eq:bracket}) 
and the following formula:

\begin{equation}
\theta \per (U\wedge V)=(U\per \theta )\per V+(-1)^{uv}(V\per \theta )
\per U .
\end{equation}
\qed

Conversely, the connection $\nabla$ can be readily 
 recovered from the operator $D$. More precisely,
we have

\begin{pro}
Suppose that  $D: \gm (\wedge^{k}A) \lon \gm (\wedge^{k-1}A)$ is
the  operator corresponding to an $A$-connection
$\nabla$ on $\wedge^{n}A$.
Then, for any $X\in \gm (A)$ and $\Lambda \in \gma{n} $,
\begin{equation}
\label{eq:connection} 
\nabla_{X}\Lambda =-X\wedge D\Lambda  .
\end{equation}
\end{pro}
\pf By definition,
$D\Lambda =-\alpha_{i}\per \nabla_{X_{i}}\Lambda$.
Hence,
\be
-X\wedge D\Lambda &=&\sum_{i} X\wedge (\alpha_{i}\per \nabla_{X_{i}}\Lambda)
\\
&=& \sum_{i}\alpha_{i}(X) \nabla_{X_{i}}\Lambda \\
&=&\nabla_{X}\Lambda  ,
\ee
where the last equality  uses the identity: $X=\sum_{i}\alpha_{i}(X)X_{i}$,
and the second equality  follows from the following simple
fact in linear algebra:

\begin{lem}
Let $V$ be any vector space, $X\in V$, $\alpha \in V^*$ and
$\Lambda \in \wedge^{n}V$. Then,
$$ X\wedge (\alpha \per \Lambda )=\alpha (X) \Lambda  .$$
\end{lem}
\qed

\begin{thm}
\label{thm:exact}
Let $A$ be a Lie algebroid with anchor $a$, and 
$\cala =\oplus_{i}\gma{i}$
its corresponding Gerstenhaber algebra.
There is a  one-to-one correspondence between
algebroid $A$-connections  on $E \cong \wedge^{n}A$ and linear 
operators $D$ generating   the Gerstenhaber
algebra bracket  on $\cala$.
Under this correspondence, flat connections 
correspond to operators of square zero.
\end{thm}
\pf It remains to prove that Equation (\ref{eq:connection})
indeed defines an $A$-connection on $\wedge^{n}A$ if $D$
is an operator generating the \Ger algebra bracket.

First, it is clear that, with this definition,
 $\nabla_{fX}\Lambda =f\nabla_{X}\Lambda$ for any $f\in C^{\infty}(M)$.

To prove that it satisfies the second axiom of a linear connection,
we observe that
for any $f\in C^{\infty}(M)$,
and $\Lambda \in \gma{n}$,
\be
D(f\Lambda )&=&(Df)\Lambda +f D\Lambda +[f, \Lambda ]\\
&=&fD\Lambda +[f, \Lambda ].
\ee
Hence,
\be
\nabla_{X}(f\Lambda ) &=&-X\wedge D(f \Lambda )\\
&=&-X\wedge (fD\Lambda +[f, \Lambda ] ) \\
&=&f\nabla_{X}\Lambda -X\wedge [f, \Lambda ] .
\ee

On the other hand, using  the property of Gerstenhaber algebras,
\be
[f, X\wedge \Lambda ]&=&[f, X]\wedge \Lambda +(-1)X\wedge [f, \Lambda]\\
&=&-(a (X)f) \Lambda -X\wedge [f, \Lambda ].
\ee
Thus, $X\wedge [f, \Lambda ]=-(a (X)f) \Lambda $. Hence,
$\nabla_{X}(f \Lambda )=f \nabla_{X}\Lambda +(a (X)f)\Lambda $.
\qed

A flat $A$-connection always exists on the line bundle $E=\wedge^{n}A$.
To see this, note that $E\otimes E$ is a trivial line bundle, which 
always admits a flat connection. So the ``square root" of this
connection (see Proposition 4.3 in \cite{EvansLW})  is
a flat connection  we need. Therefore, for a given  Lie algebroid,
there always exists an operator of degree $-1$ and of square zero
generating the corresponding \Ger algebra. Such an 
operator is called a {\em \Koszul operator}.

Any  $A$-connection $\nabla$  on  the  Lie algebroid $A$  itself
induces an $A$-connection  on the line bundle $E=\wedge^{n}A$. Therefore, 
it   corresponds to  a linear
operator  $D$ generating the Gerstenhaber algebra $\cala$.
In particular, if it is torsion free, i.e.,
$$\nabla_{X}Y-\nabla_{Y}X=[X, Y], \ \ \ \forall X, Y\in \gm (A), $$
  $D$ possesses  a simpler expression. Note that $\nabla$ induces an $A$-connection 
on the  dual bundle $A^*$, and  its
 exterior powers, which is  denoted by the same symbol.

\begin{pro}
Suppose that $\nabla$ is a torsion free $A$-connection
on $A$.  Let $D: \secc{*}\lon \secc{*-1}$ be
the  induced operator.  Then,
for any $U\in \gma{u}$,
\begin{equation}
DU|_{m}=- \sum_{i}\alpha_{i}\per \nabla_{X_{i}}U,
\end{equation}
where $X_{1}, \cdots , X_{n}$ is a basis of $A|_{m}$ and
$\alpha_{1}, \cdots , \alpha_{n}$ the  dual basis of $A^{*}|_{m}$.
\end{pro}
 \pf Assume that $U=\omega \per \Lambda$ for some $\Lambda \in \gma{n}$ and
$\omega \in \gmd{n-u}$. Then,
\be
\sum_{i}\alpha_{i}\per \nabla_{X_{i}} (\omega \per \Lambda )&=&
\sum_{i}\alpha_{i} \per [\nabla_{X_{i}}\omega \per \Lambda +
\omega \per  \nabla_{X_{i}}\Lambda ]\\
&=&\sum_{i} [ (\nabla_{X_{i}}\omega  \wedge \alpha_{i} )\per \Lambda
+(\omega \wedge \alpha_{i} )\per \nabla_{X_{i}}\Lambda ]\\
&=&(-1)^{|\omega |}(\sum_{i}(\alpha_{i} \wedge \nabla_{X_{i}}\omega ) \per
\Lambda +\sum_{i}(\alpha_{i}\wedge \omega )\per \nabla_{X_{i}}\Lambda ) .
\ee

The conclusion thus follows from the following

\begin{lem}
For any $\omega \in \gmd{|\omega |}$,
$$d\omega =\sum_{i}\alpha_{i}\wedge \nabla_{X_{i}}\omega .$$
\end{lem}
\pf Define an operator $\delta : \gmd{k}\lon \gmd{k+1}$, for all
$0\leq k\leq n$,  by
$$\delta \omega =\sum_{i}\alpha_{i}\wedge \nabla_{X_{i}}\omega .$$
It is simple to check that $\delta $ is a derivation with
respect to the wedge product, i.e., 
$$\delta (\omega \wedge \theta )=\delta \omega \wedge \theta
+(-1)^{|\omega |}\omega \wedge \delta \theta .$$

For any $f\in C^{\infty}(M)$,
$$\delta f=\sum_{i}\alpha_{i}\nabla_{X_{i}}f=\sum_{i}[a (X_{i})f]\alpha_{i}
=df .$$

For any $\theta \in \gm (A^{*}) $ and $X, \ Y\in \gm (A)$,
\be
(\delta \theta)(X, Y) &=&\sum_{i}(\alpha_{i}\wedge \nabla_{X_{i}}\theta )(X, Y)\\
&=&\sum_{i} \alpha_{i}(X)(\nabla_{X_{i}}\theta )(Y)-\alpha_{i}(Y)
(\nabla_{X_{i}}\theta )(X)\\
&=&\sum_{i}\alpha_{i}(X)(\nabla_{X_{i}}(\theta \cdot Y)-\theta \cdot
\nabla_{X_{i}}Y)-\sum_{i}\alpha_{i}(Y)(\nabla_{X_{i}}(\theta \cdot X)-
\theta \cdot \nabla_{X_{i}}X)\\
&=&\sum_{i}\alpha_{i}(X)(a (X_{i})(\theta \cdot Y)-\theta \cdot 
\nabla_{X_{i}}Y)-\sum_{i}\alpha_{i}(Y)(a (X_{i})(\theta \cdot X)- 
\theta \cdot \nabla_{X_{i}}X)\\ 
&=&a (X)(\theta \cdot Y)-\theta \cdot \nabla_{X}Y-
a (Y)(\theta \cdot X)+\theta  \cdot \nabla_{Y}X\\
&=&a (X)(\theta \cdot Y)-a (Y)(\theta \cdot X)-\theta \cdot
( \nabla_{X}Y-  \nabla_{Y}X)\\
&=&a (X)(\theta \cdot Y)-a (Y)(\theta \cdot X)-\theta \cdot [X, Y]\\
&=&d\theta (X, Y).
\ee
Therefore,   $\delta$ coincides with the exterior
derivative $d$,  since $\gma{*}$ is generated by $\gm (A^{*}) $
over the module $C^{\infty}(M)$.
\qed
{\bf Remark} (1) Theorem  \ref{thm:exact} was proved by Koszul for the case 
of  the tangent bundle  Lie algebroid $TP$ \cite{Koszul}. In fact, his   
result was the main  motivation of  our  work here.
However, Koszul used an indirect argument instead of  using
Equations (\ref{eq:D})  and 
(\ref{eq:connection}). We will see more applications of
these equations in the next section.

(2)  A flat  $A$-connection on a vector bundle $E$  is also called
a representation of the  Lie algebroid by Mackenzie \cite{Mackenzie}
  \cite{EvansLW}.

We end this section by introducing the notion of
generalized divergence. Let $\nabla$ be a flat $A$-connection
on $\wedge^{n}A$, and $D$ its corresponding \Koszul operator.
For any section $X\in \gm (A)$, we use $div_{\nabla}X$
to denote the function  $DX$. When $A=TP$ with the
usual  Lie algebroid structure, and $\nabla$ is the
flat connection induced by a volume, $DX$ is the
divergence in the ordinary  sense.  So $DX$  
 can be indeed  considered as a generalized divergence.

The following proposition gives a simple geometric
characterization for  the divergence of a section $X$ of $A$.

\begin{pro}
For any $X\in \gm (A)$ and $\Lambda \in \gm (\wedge^{n}A)$,
$$L_{X}\Lambda -\nabla_{X}\Lambda =(div_{\nabla}X)\Lambda. $$
In other words, $div_{\nabla}X$ is the function on $M$  defining
 the bundle  map $L_{X}  - \nabla_{X}  $  of 
 the  line bundle $\wedge^{n}A$.
\end{pro}
\pf Assume that $X=\omega \per \Lambda$ for some $\omega \in \gm (\wedge^{n-1}A^{*})$.
Then, 
$$DX=- (-1)^{|\omega |}(d\omega \per \Lambda+\sum_{i=1}^{n}(\alpha_{i}
\wedge \omega ) \per \nabla_{X_{i}}\Lambda ). $$

Now

\be
\sum_{i=1}^{n}( ( \alpha_{i}\wedge \omega )
\per \nabla_{X_{i}}\Lambda )\Lambda &=&
\sum_{i=1}^{n}( ( \alpha_{i}\wedge \omega ) 
\per \Lambda )\nabla_{X_{i}}\Lambda \\
&=&\sum_{i} (-1)^{|\omega |}((\omega \wedge \alpha_{i} ) \per \Lambda ) 
\nabla_{X_{i}}\Lambda\\
&=&\sum_{i} (-1)^{|\omega |}(\alpha_{i}\per X)\nabla_{X_{i}}\Lambda\\
&=&\sum_{i} (-1)^{|\omega |}X(\alpha_{i})\nabla_{X_{i}}\Lambda\\ 
&=& (-1)^{|\omega |} \nabla_{X}\Lambda.
\ee

On the other hand, it follows from Equation (\ref{eq:bracket})
 that $$[X, \Lambda ]\per \theta =-\Lambda \per d(X\per \theta ), $$
for any  $\theta \in \gmd{n}$.
It is simple to see that $X\per \theta =(\omega \per  \Lambda )\per
\theta =(-1)^{|\omega |(n- |\omega |)}\omega $. Since
$n=|\omega | -1$, then $X\per \theta= (-1)^{|\omega |} \omega$,
and $[X, \Lambda ]\per \theta =- (-1)^{|\omega |} d\omega \per \Lambda $.

Hence, $(DX  )\Lambda =[X, \Lambda ]-\nabla_{X}\Lambda 
=L_{X}\Lambda -\nabla_{X}\Lambda $. \qed

\section{Homology of Lie algebroids}
Suppose that $A$ is a Lie algebroid,
and $\nabla$  a flat $A$-connection
on the line bundle $E=\wedge^{n}A$.
Let $D$ be the corresponding
\Koszul operator  and $\partial =(-1)^{n-k}D:
\gma{k}\lon \gma{k-1}$. Then, $\partial^{2}=0$. The reason for choosing
this sign in the definition of $\partial$ will
become clear later (see Equation (\ref{eq:star})).

As usual, define the  homology by
$$H_{*}(A, \nabla )=ker \partial/Im \partial .$$

Since $D$ is a derivation   with respect to $[\cdot , \cdot ]$,
immediately we have

\begin{pro}
The Schouten bracket passes to the homology $H_{*}(A, \nabla )$. 
\end{pro}

Since this homology depends on the choice of the connection
$\nabla$, it is natural to ask  how  $H_{*}(A, \nabla )$ changes
according to the connection $\nabla$.

\begin{pro}
Let $\tilde{D}$ and $D$ be two \Koszul operators.
Then $D-\tilde{D}=i_{\alpha } $,  where $\alpha \in \gm (A^{*})$.
In this case, $$\tilde{D}^{2}-D^{2}=-i_{d\alpha }. $$
In particular, if $\tilde{D}^{2}=D^{2}=0$, then $\alpha\in  \gm (A^{*})$
is closed.
\end{pro}
\pf Let $\tilde{\nabla}$ and $\nabla $ be any two $A$-connections
on $E=\wedge^{n}A$. Then there is $\alpha \in \gm (A^{*})$ such
that
$$\tilde{\nabla} _{X} s=\nabla_{X}s+<\alpha , X>s, \ \ \forall s\in \gm (\wedge^{n} A) .  $$

Let $\tilde{D}$ and $D$ be their corresponding \Koszul operators.
It follows from a direct verification that
$$\tilde{D}=D-i_{\alpha }. $$
According to Proposition  \ref{pro:D2}, 
$\tilde{D}^{2}U-D^{2}U =-(\tilde{R}-R)\per U$,
where $\tilde{R}$ and $R$ are the curvatures of $\tilde{\nabla}$ 
and $\nabla$, respectively.
Finally, it is  routine to check that  
$\tilde{R}-R=d\alpha $.
\qed

\begin{defi}
$A$-connections $\nabla_{1}$ and  $\nabla_2$ are
said to be \equivalent if they differ by an exact form in $\gm (A^{*})$.
 
Similarly two \Koszul operators $D_{1}$ and $D_{2}$ are said
to be \equivalent if they differ by  an exact form,
i.e., $D_{1}-D_{2}=i_{\alpha }$ for some exact form $\alpha \in \gm (A^{*})$.
\end{defi}

The following result is thus immediate.
\begin{pro}
Let $\nabla_{1}$ and $\nabla_{2}$  be two flat $A$-connections
on the canonical  line bundle  $E=\wedge^{n}A $, and $D_{1}$ and $D_{2}$ 
the corresponding \Koszul operators.
If $\nabla_{1}$ and  $\nabla_2$ are \equivalent
(or equivalently  $D_{1}$ and $D_{2} $ are \equivalent), then,
\begin{equation}
  H_{*}(A, \nabla_{1} )\cong H_{*}(A, \nabla_{2} ) .
\end{equation}
\end{pro}

Now assume that $\wedge^{n} A$ is a trivial bundle,
so  there exists a nowhere vanishing volume
$\Lambda \in \gma{n}$.  This volume induces  
a flat $A$-connection  $\nabla_{0}$  on $\wedge^{n}A $ simply by
$(\nabla_{0})_{X}\Lambda =0$ for all $X\in \gm (A) $. Let $D_{0}$
be its corresponding \Koszul operator. 
Note that $\Lambda$ being  horizontal  is equivalent to the condition:
$$D_{0} \Lambda =0 . $$ 
Suppose  that $\Lambda '$ is another nonvanishing   volume, and $\nabla '$
its corresponding flat connection on $E$.
Assume that  $\Lambda '=f \Lambda$ for some  positive
$f \in C^{\infty}(M)$. Then,
it is easy to see that
$$\nabla  ' _{X} s=(\nabla_{0})_{X}s-< d\ln{f} , X>s . $$
In other words, their corresponding \Koszul operators are   \equivalent.

Let us now  fix  such a volume $\Lambda \in \gma{n}$.
Define a $*$-operator from $\gmd{k}$ to
$\gma{n-k}$ by
$$* \omega =\omega \per \Lambda  .$$
Clearly  $*$ is an isomorphism.  

The following proposition follows  immediately  from definition.

\begin{pro}
The operator  $\partial_{0}=(-1)^{n-k}D_{0}$   equals to
$-*\smalcirc d \smalcirc *^{-1}$.   That is, 
\begin{equation}
\label{eq:star}
\partial_{0}= -*\smalcirc d \smalcirc *^{-1}. 
\end{equation}
\end{pro}

Thus, as a consequence, we have

\begin{thm}
\label{thm:iso}
Let $\nabla_{0}$ be an $A$-connection on $\wedge^{n}A$ such that
  there exists a nowhere vanishing  horizontal 
volume $\Lambda \in  \gma{n}$. Then
$$H_{*}(A, \nabla_{0})\cong H^{n-*}(A, \reals ). $$
\end{thm}
{\bf Remark} We see, from the discussion above, that there is
a family of Lie algebroid homologies,  which are in
some sense  parameterized by the first Lie algebroid cohomology.
In the case that the canonical line bundle $\wedge^{n}A$
is trivial,  one of them
is isomorphic to the Lie algebroid cohomology with trivial coefficients,  and
all the others
can be considered as Lie  algebroid cohomology with some
twisted coefficients. 
In general,  the Lie algebroid homology introduced above  is
 a special case of Lie algebroid cohomology with
general coefficients  in a line bundle (see \cite{Mackenzie} \cite{EvansLW}).

Below are 
 two special interesting cases.

 Let $\frakg$ be an $n$-dimensional  Lie algebra.
Then $\wedge^{n}\frakg$ is one-dimensional and  obviously has a 
trivial $\frakg$-connection. This induces an operator
$D_{0}: \wedge^{*}\frakg \lon \wedge^{*-1}\frakg$  which
is of square zero and generates the Schouten bracket on $\wedge^{*} \frakg $.
On the other hand, there exists another operator 
$D:\wedge^{*}\frakg \lon \wedge^{*-1}\frakg$, which is dual to the differential
 operator of the Lie algebra cohomology.
In general, $D$ is different from $D_{0}$ and  in fact 
it is easy to check  that $D-D_{0}=i_{\alpha}$, 
where $\alpha $ is the modular character of the Lie algebra.
In particular, when $\frakg$ is a unimodular Lie algebra,
the Lie algebra homology is isomorphic to Lie algebra
cohomology, a well-known result.

Another interesting case, which does not seem trivial, is
the one when $A$ is the cotangent Lie algebroid $T^{*}P$ 
of a Poisson manifold $P$ (see Equations (\ref{eq:anchor0}) and
 (\ref{eq_bracket-on-one-forms})). 
In this case, $\gm (\wedge^{k}T^{*}P) =\Omega^{k}(P)$.
There exists  an  operator $D: \Omega^{k}(P)\lon \Omega^{k-1}(P)$
 introduced by  Koszul \cite{Koszul}  and studied by Brylinski \cite{Br}.
 It is given by
$$D=[i_{\pi} , d].  $$

The corresponding homology is called Poisson homology by
Brylinski, and denoted by $H_{*}(P, \pi )$.

It was shown in \cite{Koszul} that the
operator $D$  indeed  generates the \Ger bracket on  
$\Omega^{*}(P)$ induced from  the cotangent Lie algebroid of $P$.

Therefore, $D$  corresponds to a flat Lie algebroid  connection
on $\wedge^{n}T^{*}P$, which is  given by

\begin{equation}
\label{eq:connection1}
\nabla_{\theta }\Omega =-\theta \wedge D\Omega=\theta \wedge d(\pi \per \Omega ), 
\end{equation}
for any $\theta \in \Omega^{1}(P) $ and $\Omega \in \Omega^{n}(P)$, according
to Equation (\ref{eq:connection}).  A similar formula
  was also discovered independently,
by Evens-Lu-Weinstein \cite{EvansLW}.

The  Koszul-Brylinski operator $D$ is intimately related to
the so called modular class of the Poisson manifold, a classical
analogue of  the modular form of a von Neumann  algebra, which was 
introduced  recently by Weinstein \cite{we:modular}, and
independently by Brylinski and  Zuckerman \cite{BZ}.

For simplicity, let us assume that $P$ is orientable,
and $\Omega $ is a   volume form.  The modular vector
field $\nu_{\Omega}$ is the vector field defined  by
$$f\lon (L_{X_{f}}\Omega )/\Omega,  \ \ \ \ \forall f\in C^{\infty}(P).$$
It can be shown that the above map is a derivation
on the space of functions $C^{\infty}(P)$, so it indeed 
defines a vector field. It also can be shown that 
$\nu_{\Omega}$ is a Poisson vector field.  When the volume $\Omega$ changes, the
corresponding modular vector fields differ  by   a hamiltonian
vector field. So its class is a well-defined element in
the first Poisson cohomology $H^{1}_{\pi }(P)$, which
is called the modular class of the Poisson manifold.
A Poisson manifold is called {\em unimodular} if
its modular class vanishes.
In fact, the modular class can be defined for 
any Poisson manifold by just  replacing the volume by a positive
density. We refer the  interested reader to \cite{we:modular}
for more detail.

Now let $P$ be  an  orientable Poisson manifold with 
  volume form $\Omega$, and let $D_{0}$ be its corresponding \Koszul operator
as in the observation preceding Theorem \ref{thm:iso}.

The following proposition follows  immediately
from a direct verification.

\begin{pro}
\label{pro:modular}
Let $D$ be the Koszul-Brylinski operator of a Poisson manifold $P$.
Then
 $D-D_{0}=i_{\nu_{\Omega}}$, where $\nu_{\Omega}$ is the modular vector field
corresponding to the volume $\Omega$.
\end{pro} 

As an immediate consequence, we have

\begin{thm}
If   $P$ is an orientable unimodular  Poisson manifold, then 
$$H_{*}(P, \pi )\cong H^{n-*}_{\pi} (P). $$
\end{thm}

In particular, this holds for any symplectic manifold, which 
 was first   proved by Brylinski \cite{Br}. \\\\\\
{\bf Remark}  The above situation can be generalized to
the case of triangular Lie bialgebroids. 
Let $A$ be a Lie algebroid with anchor $a$.   A triangular $r$-matrix
is a section $\pi $ in $\gma{2}$ satisfying the condition
$[\pi, \pi ]=0$. One may think that this is a sort of
generalized ``Poisson structure" on the generalized
manifold $A$. In this case, $A^*$ is equipped with 
a Lie algebroid structure with the anchor $a\smalcirc \pi^{\#}$
and the Lie bracket defined by an equation  identical to
 the one defining
the bracket on one-forms of a Poisson manifold.

Similarly, $D=[i_{\pi}, d]$ is an operator  $\gmd{k}\lon \gmd{k-1}$
of square zero and generates the bracket $[\cdot , \cdot ]$
on $\gmd{*}$.
A form of top degree $\Omega \in \gmd{n}$ satisfies
the condition $D\Omega =0$ iff $\pi \per \Omega \in \gmd{n-2}$
is closed. If there exists such a nowhere vanishing  
form, the homology $H_{*}(A, \nabla  )$  is then 
isomorphic  to the cohomology $H^{n-*}(A, \reals )$.

\section{Discussion}

We end this paper  by some  open questions.

{\bf Question 1}: As in the remark at the end
of the last section,  is the condition that $\pi \per \Omega \in \gmd{n-2}$
is closed equivalent to that  the Lie algebroid $A^*$  is  unimodular?

{\bf Question 2}: For a general Lie algebroid $A$, does  there exist
any canonical  \Koszul operator
  corresponding  to the
modular class of  the Lie algebroid  just as in the case
of cotangent algebroid of a Poisson manifold (see  Proposition  \ref{pro:modular})?

{\bf Question  3}: For a Poisson manifold, there
is a family of homologies parameterized by
the first Poisson cohomology $H^{1}_{\pi }(P)$.
What is  the meaning of the  homologies  other than
the Poisson  homology?

{\bf Question 4}:  Suppose that $(A, A^{*})$ is a Lie bialgebroid
and $\nabla$ a  flat $A$-connection on $\wedge^{n}A$. Then $(\gma{*}, \wedge , 
d_{*}, [ , ], D )$ is a strong differential   BV-algebra. It is clear that
$d_{*}D+Dd_{*} $ is a derivation with respect to both 
$\wedge $ and $[ , ]$. When is  $d_{*}D+Dd_{*} $  inner
and in particular, when is $d_{*}D+Dd_{*} =0$?

When $A=T^{*}P$ is the cotangent Lie algebroid of a Poisson
manifold, $A^* =TP$ the usual Lie algebroid
on the tangent bundle,  and the connection $\nabla$
 is as in Equation (\ref{eq:connection1}),   then $d_{*}$ is the usual
de-Rham differential and $D$ is the Koszul-Brylinski operator.
Thus, $d_{*}D+Dd_{*} $  is automatically
zero,  which gives rise to the  Brylinski double complex.
On the other hand, if we switch the order and
 consider $A=TP$ and $A^{*}=T^{*}P$
for a Poisson manifold $P$ equipped with a volume generating
the connection on the line bundle $\wedge^{n}TP$, 
 then $\cala=\oplus_{i}\gma{i}$ is the space of  multivector fields.
In this case, $d_{*}=[\pi , \cdot ]$ is the Lichnerowicz
differential defining the Poisson cohomology, and $D=
 -(-1)^{n-k} *\smalcirc d \smalcirc *^{-1}$, where
$*$ is the isomorphism between multivector fields and
differential forms induced by the volume element.
Then   $ d_{*}D+Dd_{*}=L_{X}$, where $X$ is the modular vector 
field of the Poisson manifold (see P. 265 of \cite{Koszul}).
It vanishes  iff  $P$ is unimodular.
It would be interesting to explore this in general.

\end{document}